\documentclass[11pt]{article}

\usepackage{amsfonts,mathrsfs}
\usepackage{amssymb,amsmath}
\usepackage{url,color}
\usepackage{makeidx}         
\usepackage{graphicx}        

\usepackage{epstopdf}

\usepackage{enumerate}

\newcommand{\mW}{\mathbf{W}}

\newcommand{\oz}{\overline{z}}

\newcommand{\bE}{\mathbb{E}}

\newcommand{\T}{\intercal}

\usepackage[left=1.25in, right=1.25in, top=1in, bottom=1in]{geometry}

\usepackage{epsfig}  

{\bf}{}
{\bf}{}
{\bf}{}
\newtheorem{lemma}{Lemma}{\bf}{}
\newtheorem{theorem}{Theorem}{\bf}{}
{\bf}{}
\newtheorem{assumption}{Assumption}{\bf}{}

\usepackage{subcaption}
\usepackage{epsfig}  
\usepackage{wrapfig}

\usepackage{float}

\newcommand{\degree}{\mathrm{deg}}

\usepackage{setspace}
\title{Asymptotic Network Independence in \\ Distributed Stochastic Optimization for Machine Learning~\footnote{The research was partially supported by the NSF under grants ECCS-1933027, IIS-1914792,
		DMS-1664644, and CNS-1645681, by the
		ONR under grant N00014-19-1-2571, by the NIH under grant 1R01GM135930, and by the SRIBD Startup Fund JCYJ-SP2019090001. (Corresponding author: Shi Pu.)}}

	\author{Shi Pu\thanks{Institute for Data and Decision Analytics, The Chinese University of Hong Kong, Shenzhen, China and Shenzhen Research Institute of Big Data
			({pushi@cuhk.edu.cn}). The research was conducted when the author was with Division of Systems Engineering, Boston
			University, Boston, MA.}
		\and
		Alex Olshevsky\thanks{Department of Electrical and Computer Engineering and Division of Systems Engineering, Boston
			University, Boston, MA
			({alexols@bu.edu}, {yannisp@bu.edu}).}
		\and Ioannis Ch. Paschalidis\footnotemark[3]
	}

\begin{document}
	\date{}
	\maketitle
	
	\begin{abstract}
		We provide a discussion of several recent results which, in certain scenarios, are able to overcome a
		barrier in distributed stochastic optimization for machine learning. Our focus is the so-called asymptotic network independence property, which is achieved whenever a distributed method executed over a network of $n$ nodes asymptotically converges to the optimal solution at a comparable rate to a centralized method with the same computational power as the entire network. 
		We explain this property through an example involving the training of ML models and sketch a short mathematical analysis for comparing the performance of distributed stochastic gradient descent (DSGD) with centralized stochastic gradient decent (SGD).
	\end{abstract}
	
	\section{Introduction:  Distributed Optimization and Its Limitations} 
	
	First-order optimization methods, ranging from vanilla gradient descent to Nesterov
	acceleration and its many variants, have emerged over the past decade as the
	principal way to train Machine Learning (ML) models. There is a great need for
	techniques to train such models quickly and reliably in a distributed fashion over networks
	where the individual processors or GPUs may be scattered across the globe and
	communicate over an unreliable network which may suffer from message losses, delays, and
	asynchrony (see \cite{brisimi2018federated,scaman2019optimal,assran2018stochastic,ying2018supervised}).
	
	Unfortunately, what often happens is that the gains from having many different
	processors running an optimization algorithm are  squandered by the cost
	of coordination, shared memory, message losses and latency. This effect is especially pronounced
	when there are many processors and they are spread
	across geographically distributed data centers. As is widely recognized by the
	distributed systems community,  ``throwing'' more processors at a problem will not,
	after a certain point, result in better performance.
	
\begin{figure}
		\centering
		\includegraphics[width=0.8\textwidth]{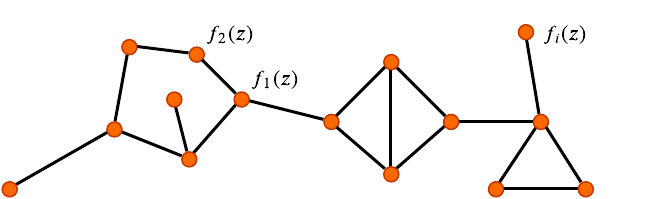}
		\caption{Example of a network. Two nodes are connected if there is an edge between them.} \label{Network}
	\end{figure}

	This is typically reflected in the convergence time bounds obtained for distributed optimization in the literature. The problem formulation is that one must solve 
	\begin{equation}
		z^* \in \arg \min_{z \in \mathbb{R}^d}\sum_{i=1}^n f_i(z),   \label{opt Problem_def}
	\end{equation}  
	over a network of $n$ nodes (see Figure \ref{Network} for an example).
	Only node $i$ has knowledge of the function $f_i(z)$,
	and the standard assumption is that, at every step when it is awake, node $i$ can
	compute the (stochastic) gradient of its own local function $f_i(z)$. These functions $f_i(z)$ are assumed to be convex. The problem is to compute this minimum in a
	distributed manner over the network based on peer-to-peer communication, possible message
	losses, delays, and asynchrony.
	
	
	This relatively simple formulation captures a large variety of learning
	problems. Suppose each agent $i$ stores training data points $\mathcal{X}_i=\{(x_j,y_j)\}$,
	where $x_j\in \mathbb{R}^p$ are vectors of features and $y_j\in \mathbb{R}$ are the
	associated responses (either discrete or continuous). We are interested to learn a
	predictive model $h(x; \theta)$, parameterized by parameters $\theta\in \mathbb{R}^d$, so that $h(x_j;
	\theta)\approx y_j$ for all $j$. In other words, we are looking for a model that fits all the data throughout the network. This can be accomplished by empirical risk minimization
	\begin{equation}
		\label{Empirical_rm}
		\theta^* \in \arg\min_{\theta \in \mathbb{R}^d} \sum_{i=1}^n c_i(\theta,\mathcal{X}_i),
	\end{equation}
	where 
	\[c_i(\theta,\mathcal{X}_i)=\sum_{(x_j,y_j)\in\mathcal{X}_i}\ell(h(x_j; \theta), y_j)\]
	measures how well the parameter $\theta$ fits the data at node $i$,
	with $\ell(h(x_j; \theta), y_j)$ being a loss function measuring the difference between $h(x_j;
	\theta)$ and $y_j$.
	Much of modern machine learning is built around such a
	formulation, including regression, classification, and regularized variants \cite{chen2018robust}.
	
	It is also possible that each agent $i$ does not have a static dataset,
	but instead collects streaming data points $(x_i,y_i)\sim \mathbb{P}_i$ repetitively over time, where $\mathbb{P}_i$ represents an unknown distribution of $(x_i,y_i)$. In this case we can find $\theta^*$ through expected risk minimization
	\begin{equation}
		\label{Expected_rm}
		\theta^* \in \arg\min_{\theta \in \mathbb{R}^d}\sum_{i=1}^n f_i(\theta),
	\end{equation}
	where
	\[f_i(\theta)=\bE_{(x_i,y_i)\sim \mathbb{P}_i}\ell(h(x_i; \theta), y_i).\]
	
	This paper is concerned with the current limitations of  distributed optimization and how to get past them in certain scenarios. To illustrate our main concern, let us consider the distributed subgradient method in the simplest possible setting, namely the problem of computing the median of a collection of numbers in a distributed manner over a fixed graph.  Each agent $i$ in the network holds value $m_i>0$, and the global objective is to find the median of $m_1,m_2,\ldots,m_n$.   This can be incorporated in the framework of (\ref{opt Problem_def}) by choosing 
	\[f_i(z)=|z-m_i|, \quad \forall i.\] 
	The distributed subgradient method (see \cite{nedic2009distributed2}) uses the subgradients $s_i(z)$ of $f_i(z)$ at any point $z$,  to have agent $i$ update as 
	\begin{equation}
		\label{eq: z_i,k_subgradient}
		z_i(k+1) = \sum_{j=1}^{n}w_{ij}z_j(k)-\alpha_k s_i(z_i(k)),
	\end{equation}
	where $\alpha_k>0$ denotes the stepsize at iteration $k$, and $w_{ij}\in [0,1]$ are the weights agent $i$ assigns to agent $j$'s solutions: two agents $i$ and $j$  are able to exchange information if and only if $w_{ij}, w_{ji}>0$ ($w_{ij}=w_{ji}=0$ otherwise). The weights $w_{ij}$ are assumed to be symmetric. For comparison, the centralized subgradient method updates the solution at iteration $k$ according to
	\begin{equation}
		\label{eq: z_k_subgradient}
		z(k+1) = z(k)-\alpha_k \frac{1}{n}\sum_{j=1}^{n}s_j(z(k)).
	\end{equation}
	
\begin{figure}[htbp]
		\centering
		\includegraphics[width=0.6\textwidth]{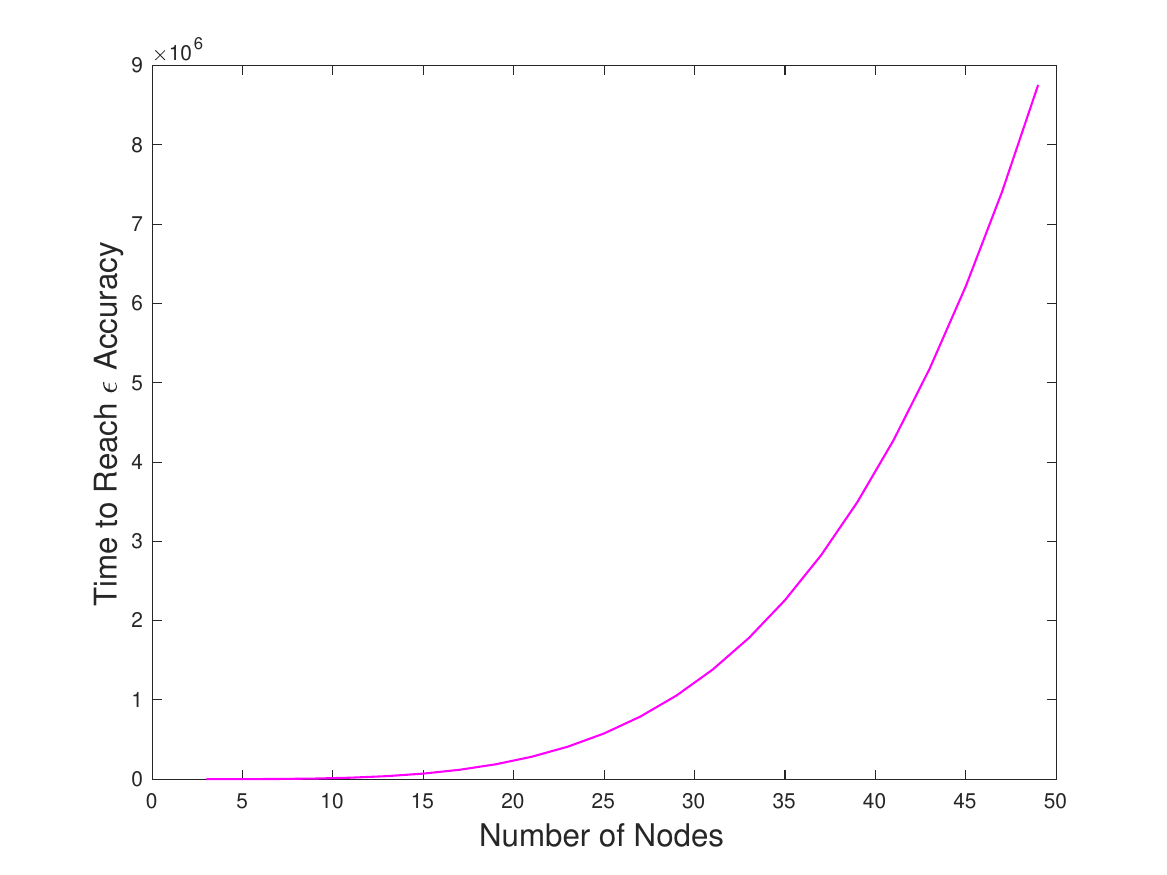}
		\caption{Performance of Algorithm (\ref{eq: z_i,k_subgradient}) as a function of the network size $n$. The agents communicate over a ring network (see Figure \ref{Ring}) and choose the Metropolis weights (see Section \ref{subsec: setup} for the definition). Stepsizes $\alpha_k=\frac{1}{\sqrt{k}}$, and $m_i$ are evenly distributed in $[-10,10]$. The time $k$ to reach $\frac{1}{n}\sum_{i=1}^n|y_i(k)|<\epsilon$ is plotted, where $y_i(k)=\frac{1}{k}\sum_{\ell=0}^{k-1} z_i(\ell)$ and $\epsilon=0.1$.} \label{Median}
	\end{figure}

	In Figure \ref{Median}, we show the performance of Algorithm (\ref{eq: z_i,k_subgradient}) as a function of the network size $n$ assuming the agents communicate over a ring network. As can be clearly seen, when the  network size grows it takes a longer time for the algorithm to reach a certain performance threshold.
	
	Clearly, this is an undesirable property. Glancing at the figure, we see that distributing computation over $50$ nodes can result in a convergence time on the order of $10^7$ iterations. Few practitioners will be enthusiastic about distributed optimization if the final effect is vastly increased convergence time. 
	
	One might hope that this phenomenon, demonstrated for the problem of median computation -- considered here because it is arguably the simplest problem to which one can apply the subgradient method -- will not hold for the more sophisticated optimization problems in the ML literature. Unfortunately, most work in distributed optimization replicates this undesirable phenomenon. We next give an extremely brief discussion of known convergence times in the distributed setting (for a much more extended discussion, we refer the reader to the recent survey \cite{nor}). 
	
	We would like to confine our discussion to the following point: most known convergence times in the distributed optimization literature imply bounds of the form 
	\begin{equation} \label{eq:pdef} 
		{\rm Time}_{n,\epsilon}({\rm decentralized})  \leq  p(\mathcal{G}) {\rm Time}_{n,\epsilon}({\rm
			centralized)}, 
	\end{equation}
	where ${\rm Time}_{n,\epsilon}({\rm decentralized})$ denotes the time for the
	decentralized algorithm on $n$ nodes to reach $\epsilon$ accuracy (error $<\epsilon$), and ${\rm Time}_{n,\epsilon}({\rm centralized)}$ is the time for the centralized algorithm
	{\em which can query $n$ gradients per time step} to reach the same level of accuracy.  The graph $\mathcal{G}=(\mathcal{N},\mathcal{E})$ consists of the set of nodes and edges in the network, denoted by $\mathcal{N}$ and $\mathcal{E}$, respectively. The function $p(\mathcal{G})$ can usually be bounded in terms of some polynomial in the number of nodes $n$. 
	
	For instance, in the subgradient methods, Corollary 9 of \cite{nor} implies that
	\begin{align*}
		& {\rm Time}_{n,\epsilon}({\rm decentralized})=\mathcal{O}\left(\frac{\max\{\|\frac{1}{n}\sum_{i=1}^n z_i(0)-z^*\|^2,G^4 h(\mathcal{G})\}}{\epsilon^2}\right),\\
		& {\rm Time}_{n,\epsilon}({\rm centralized)}=\mathcal{O}\left(\frac{\max\{\|z(0)-z^*\|^2,G^4\}}{\epsilon^2}\right), 
	\end{align*}
	where $z(0),z_i(0)$ are initial estimates, $z^*$ denotes the optimal solution and $G$ bounds the $\ell_2$-norm of the subgradients. The function $ h(\mathcal{G})$ is the inverse of the spectral gap corresponding to the graph, and will typically grow with $n$; hence, when $n$ is large, $p(\mathcal{G})\simeq h(\mathcal{G})$. In particular
	if the communication graphs are 1) path graphs, then $p(\mathcal{G})=\mathcal{O}(n^2)$;
	2) star graphs, then $p(\mathcal{G}) =\mathcal{O}(n^2)$;
	3) geometric random graphs, then $p(\mathcal{G})=\mathcal{O}(n\log n)$.
	The method developed in \cite{olshevsky2017linear} achieves $p(\mathcal{G}) = n$, but typically $p(\mathcal{G})$ is at least $n^2$.
	
	By comparing ${\rm Time}_{n,\epsilon}({\rm decentralized}) $ and ${\rm Time}_{n,\epsilon}({\rm
		centralized)}$, we are keeping the computational power the same in both
	cases. Naturally, the centralized is always better: anything that can be done in a
	decentralized way could be done in a centralized way. The question, though, is {\em
		how much better.}
	
	Framed in this way, the polynomial scaling in the quantity $p(\mathcal{G})$ is extremely
	disconcerting. It is hard, for example, to argue that an algorithm should be run in a
	distributed manner with, say, $n=100$ if  the quantity $p(\mathcal{G})$ in Eq. (\ref{eq:pdef}) satisfies $p(\mathcal{G})=n^2$; that would imply the
	distributed variant would be $10,000$ times slower than the centralized one with the same computational power.
	
	Sometimes $p(\mathcal{G})$ is written as the inverse spectral gap $\frac{1}{1-\lambda_2}$ in terms of the
	second-eigenvalue of some matrix. Because the second-smallest eigenvalue of an undirected graph Laplacian
	is approximately $\sim 1/n^2$ away from zero, such bounds will translate into at least
	quadratic scalings with $n$ in the worst-case. Over time-varying $B$-connected graphs, the best-known
	bounds on $p(\mathcal{G})$ will be cubic in $n$ using the results of
	\cite{nedic2009distributed}.
	
	There are a number of caveats to the pessimistic argument outlined above. For
	example, in a multi-agent scenario where data sharing is not desirable or feasible,
	decentralized computation might be the only available option. Generally speaking, however,
	fast-growing $p(\mathcal{G})$ will preclude the widespread applicability of distributed
	optimization. Indeed, returning to the back-of-the-envelope calculation above, if a
	user has to pay a {\em multiplicative factor of 10,000} in convergence speed to use
	an algorithm, the most likely scenario is that the algorithm will not be used.
	
	There are some scenarios which avoid the pessimistic discussion above: for example, when the underlying graph is an 
	expander, the associated spectral gap is constant (see Chapter 6 of \cite{durrett2007random} for a definition of these terms as well as an explanation), and likewise when the graph is a star graph. In particular, 
	on a random Erd\H{o}s-R\'{e}nyi random graph, the quantity $p(\mathcal{G})$ is constant with high probability (Corollary 9, part 9 in \cite{nor}). Unfortunately, these are very special cases which may not always be realistic. A star graph requires a single node to have the ability to receive and broadcast messages to all other nodes in the system. On the other hand, an expander graph may not  occur in geographically distributed systems. By way of comparison, a random graph where nodes
	are associated with random locations, with links between nodes close together, will not have constant spectral gap and will thus have $p(\mathcal{G})$ that grows with $n$
	(Corollary 9, part 10 of \cite{nor}). The Erd\H{o}s-R\'{e}nyi graph escapes this because, if we again associate nodes with locations, the average link in such a graph
	is a ``long range'' one connecting nodes that are geographically far apart.  It is a consequence of Cheeger's inequality that graphs based on connecting nearest neighbors (i.e., where nodes are regularly spaced in $\mathbb{R}^d$ and each node is connected to a constant number of closest neghbors) will not have
		constant spectral gap.
	

	\section{Asymptotic Network Independence in Distributed Stochastic Optimization} 
	
	In this paper, we provide a discussion of several recent papers which have obtained that,
	for a number of settings involving distributed stochastic optimization, $p(\mathcal{G})=1$, as long as $k$ is large enough. In other words,
	asymptotically, the distributed stochastic gradient algorithm converges to the optimal solution at a comparable rate to a centralized algorithm
	with the same computational power.

	We call this property {\em asymptotic network independence}: it is as if the network
	is not even there. Asymptotic network independence provides an answer to the
	concerns raised in the previous section.
	
\begin{figure}
		\centering
		\begin{subfigure}[t]{0.4\textwidth}
			\includegraphics[width=\textwidth]{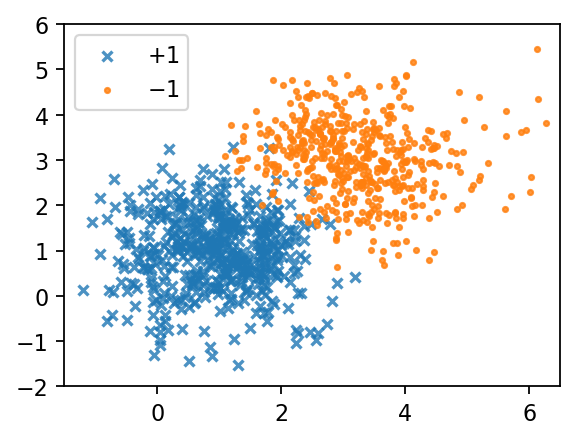}
			\caption{A total number of $1000$ data points and their labels for SVM classification. The data points are randomly generated around $50$ cluster centers.} \label{A}
		\end{subfigure}
		\begin{subfigure}[t]{0.45\textwidth}
			\includegraphics[width=\textwidth]{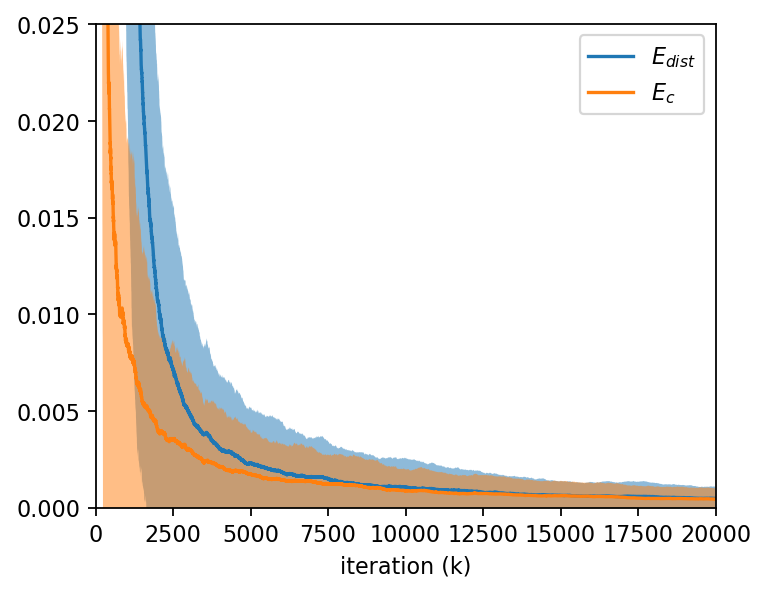}
			\caption{Squared errors and $1$-standard-deviation bands for distributed and centralized stochastic
				gradient descent. The
				performance of the centralized algorithm is shown in orange, and the performance of the
				decentralized algorithm is shown in dark blue. A total of $1000$ Monte-Carlo simulations are performed for estimating the average performance.} \label{B}
		\end{subfigure}
		\caption{Comparison between distributed and centralized stochastic gradient descent for training an SVM.}\label{SVM}
	\end{figure}

	We begin by illustrating these results with a simulation from \cite{olshevsky2018robust}, shown
	in Figure~\ref{SVM}. Here the problem to be solved is classification with
	a smooth support vector machine between  overlapping clusters of points.  The
	performance of the centralized algorithm is shown in orange, and the performance of the
	decentralized algorithm is shown in dark blue. The graph is a ring of 50 nodes, and the
	problem being solved is the search for a support vector classifier.  
	{\em The graph
		illustrates the main result, which is that a network of 50 nodes performs as well
		in the limit as a centralized method with 50 times the computational power of one node.}
	Indeed, after $\sim 8,000$ iterations the orange and dark blue lines are almost
	indistinguishable.
	
	We mention that similar simulations are available for other machine learning methods
	(training neural networks, logistic regression, elastic net regression, etc.). The asymptotic network independence property enables us to efficiently distribute the
	training process for a variety of existing learning methods.
	
	The name ``asymptotic network independence'' is a slight misnomer, as we actually do not care if the 
	asymptotic performance depends in some complicated way on the network. All we want is that the decentralized convergence rate
	can be bounded by $O(1)$ times the convergence rate of the centralized method. 
	
	The papers \cite{chen2012limiting,chen2015learning,chen2015learning2,towfic2016excess} gave
		the first crisp statement of the relationship between centralized and distributed
		methods in the setting of distributed optimization of smooth strongly convex
		functions in the presence of noise. Under constant stepsizes, the papers \cite{chen2012limiting, chen2015learning,chen2015learning2} were the first to show that when the stepsize is sufficiently small, a distributed stochastic gradient method achieves comparable performance to a centralized method in terms of the steady-state mean-square-error. The stepsize has to be small enough as a function of the network topology for this to hold.
		The paper \cite{towfic2016excess} showed that the distributed stochastic gradient algorithm asymptotically achieves comparable convergence rate to a centralized method, but assuming that all the local functions $f_i$ have the same minimum. This gives the first ``asymptotic network independence'' result.\footnote{In this survey, all the mentioned algorithms that enjoy the asymptotic network independence property assume smooth objective functions, i.e., functions with Lipschitz continuous gradients.} 
		
	The work in \cite{pu2017flocking} approximated distributed stochastic gradient descent by stochastic differential equations in continuous time by assuming sufficiently small constant stepsize. It was shown that the distributed method outperforms a centralized scheme with synchronization overhead. However, it did not lead to straightforward algorithmic bounds. 
	In our recent work \cite{olshevsky2018robust}, we generalized the results to
	graphs which are time-varying, with delays, message losses, and asynchrony. 
	In a parallel recent work \cite{koloskova2019decentralized}, a similar result was demonstrated with a further compression technique which allowed nodes to save on communication.
	
	When the objective functions are not assumed to be convex, several recent works have obtained asymptotic network independence for distributed stochastic gradient descent. In \cite{morral2014success,morral2017success},  a general
		stochastic approximation setting was considered with decaying step-sizes, and the convergence rates of centralized and distributed methods were shown to be asymptotically the same; the proof proceeded based on certain technical properties of stochastic approximation
		methods. The work in \cite{lian2017can} was the first to show that distributed algorithms could achieve a speedup like a centralized method when the number of computing steps is large enough. Such a result was generalized to the setting of directed communication networks in \cite{assran2018stochastic} for training deep neural networks, where the push-sum technique was combined with the standard distributed stochastic gradient scheme.
	
	In the rest of this section, we will give a simple and readable explanation of
	the asymptotic network independence phenomenon in the context of distributed stochastic optimization over smooth and strongly convex objective functions.
	\footnote{For more references on the topic of distributed stochastic optimization, the readers may refer to \cite{nedic2016stochastic,lan2017communication,sayin2017stochastic,sirb2018decentralized,pu2018distributed,pu2018swarming,xin2019variance} and the references therein.}
	
	\subsection{Setup}
	\label{subsec: setup}
	
	We are interested in minimizing Eq. (\ref{opt Problem_def}) over a network of $n$ communicating agents. 
	Regarding the objective functions $f_i$ we make the following standing assumption.
	\begin{assumption}
		\label{asp: mu-L_convexity}
		Each $f_i:\mathbb{R}^d\rightarrow \mathbb{R}$ is  $\mu$-strongly convex with $L$-Lipschitz continuous gradients, i.e., for any $z,z'\in\mathbb{R}^d$,
		\begin{equation}
			\langle \nabla f_i(z)-\nabla f_i(z'),z-z'\rangle\ge \mu\|z-z'\|^2,\quad
			\|\nabla f_i(z)-\nabla f_i(z')\|\le L \|z-z'\|.
		\end{equation}
	\end{assumption}
	Under Assumption \ref{asp: mu-L_convexity}, Problem (\ref{opt Problem_def}) has a unique optimal solution $z^*$, and the function $f(z)$ defined as
	\[f(z)=\frac{1}{n}\sum_{i=1}^n f_i(z)\] 
	has the following contraction property (see \cite{qu2017harnessing} Lemma 10).
	\begin{lemma}
		\label{lem: contraction_mu-L_convexity}
		For any $z\in\mathbb{R}^d$ and $\alpha\in(0,1/L)$, we have $\|z-\alpha\nabla f(z)-z^*\|\le (1-\alpha\mu)\|z-z^*\|$.
	\end{lemma}
	In other words, gradient descent with a small stepsize reduces the distance between the current solution and $z^*$.
	
	In the stochastic optimization setting, we assume that at each iteration $k$ of the algorithm, $z_i(k)$ being
	the input for agent $i$, each agent is able to obtain noisy gradient estimates $g_i(z_i(k),\xi_i(k))$ that satisfy the following condition.
	\begin{assumption}
		\label{asp: gradient samples}
		For all $i\in\{1,2,\ldots,n\}$ and $k\ge 1$, 
		each random vector $\xi_i(k)\in\mathbb{R}^m$ is independent, and
		\begin{equation}
			\label{condition: gradient samples}
			\begin{split}
				& \mathbb{E}_{\xi_{i,k}}[g_i(z_i(k),\xi_i(k))\mid z_k] =  \nabla f_i(z_i(k)),\\
				& \mathbb{E}_{\xi_{i,k}}[\|g_i(z_i(k),\xi_i(k))-\nabla f_i(z_i(k))\|^2\mid z_i(k)]  \le  \sigma^2,\quad\hbox{\ for some $\sigma>0$}.
			\end{split}
		\end{equation}
	\end{assumption}
	Stochastic gradients appear, for instance, when the gradient estimation of $c_i(\theta,\mathcal{X}_i)$ in empirical risk minimization (\ref{Empirical_rm}) introduces noise from various sources, such as sampling and quantization errors. For another example, when minimizing the expected risk in (\ref{Expected_rm}), where  independent data points $(x_i,y_i)$ are gathered over time, $g_i(z,(x_i,y_i))=\nabla_z \ell(h(x_i; z), y_i)$ is a stochastic, unbiased estimator of $\nabla f_i(z)$ satisfying the first condition in (\ref{condition: gradient samples}). The second condition holds for popular problems such as smooth Support Vector Machines, logistic regression and softmax regression assuming the domain of $(x_i,y_i)$ is bounded.
	
	The algorithm we discuss is the {\em Distributed Stochastic Gradient Descent} (DSGD) method adapted from DGD and the diffusion strategy \cite{chen2012diffusion}; note that in \cite{chen2012diffusion} this method was called ``Adapt-then-Combine''. 
	We let each agent $i$ in the network hold a local copy of the decision vector denoted by $z_i\in\mathbb{R}^d$, and its value at iteration/time $k$ is written as $z_i(k)$.  Denote $g_i(k)=g_i(z_i(k),\xi_i(k))$ for short. At each step $k\ge 0$, 
	every agent $i$ performs the following update:
	\begin{equation}
		\label{eq: z_i,k}
		z_i(k+1) = \sum_{j=1}^{n}w_{ij}\left(z_j(k)-\alpha_k g_j(k)\right),
	\end{equation}
	where $\{\alpha_k\}$ is a sequence of nonnegative non-increasing stepsizes. The initial vectors $z_{i}(0)$ are arbitrary for all~$i$, and $\mW=[w_{ij}]$ is a mixing matrix. 
	
	DSGD belongs to the class of so-called consensus-based distributed optimization methods, where different agents mix their estimates at each iteration to reach a consensus of the solutions, i.e., $z_i(k) \approx z_j(k)$ for all $i$ and $j$ in the long run.
	To achieve consensus, the following condition is assumed on the mixing matrix and the communication topology among agents.
	\begin{assumption}
		\label{asp: network}
		The graph $\mathcal{G}$ of agents is undirected and connected (there exists a path between any two agents). The mixing matrix $\mW$ is nonnegative, symmetric and doubly stochastic, 
		i.e., $\mW\mathbf{1}=\mathbf{1}$ and $\mathbf{1}^{\T}\mW=\mathbf{1}^{\T}$, where $\mathbf{1}$ is the all one vector.
		In addition, $w_{ii}>0$ for some $i\in\{1,2,\ldots,n\}$.
	\end{assumption}
	Some examples of undirected connected graphs are presented in Figure \ref{Connected_graph} below.
	\begin{figure}[H]
		\centering
		\begin{subfigure}{0.24\textwidth}
			\includegraphics[width=\textwidth]{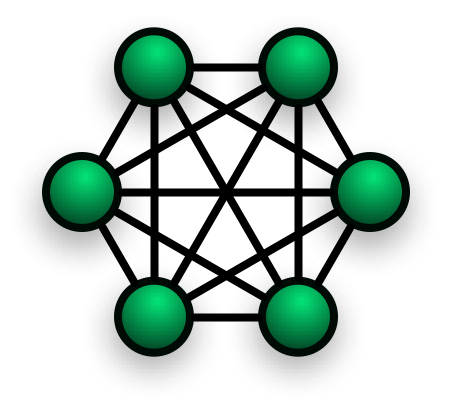}
			\caption{Fully connected graph.} \label{FullyConnected}
		\end{subfigure}
		\begin{subfigure}{0.24\textwidth}
			\includegraphics[width=\textwidth]{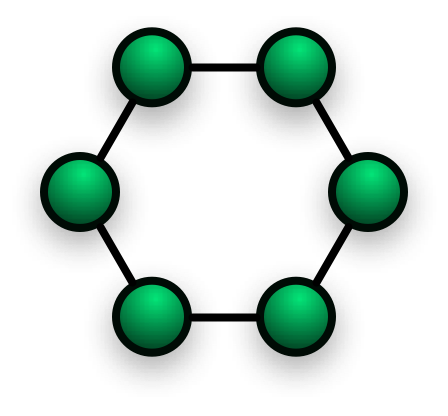}
			\caption{Ring network.} \label{Ring}
		\end{subfigure}
		\begin{subfigure}{0.24\textwidth}
			\includegraphics[width=\textwidth]{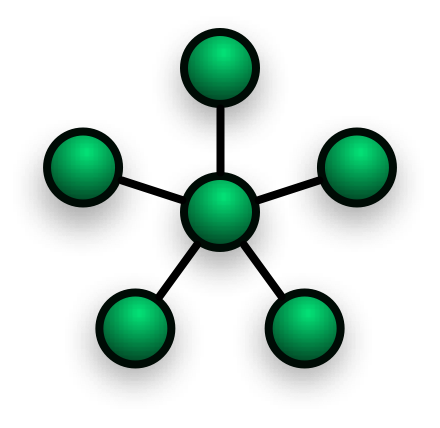}
			\caption{Star network.} \label{Star}
		\end{subfigure}
		\begin{subfigure}{0.24\textwidth}
			\includegraphics[width=\textwidth]{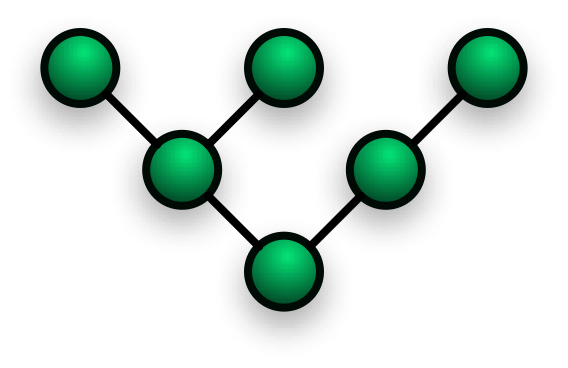}
			\caption{Tree network.} \label{Tree}
		\end{subfigure}
		\caption{Examples of undirected connected graphs.}
		\label{Connected_graph}
	\end{figure}
	Because of  Assumption \ref{asp: network}, the mixing matrix $\mW$ has an important contraction property.
	\begin{lemma}
		\label{lem: spectral norm}
		Let Assumption \ref{asp: network} hold, and let $1 = \lambda_1 \geq \lambda_2 \geq \cdots \lambda_n$ denote the  eigenvalues of 
		the matrix $\mW$. Then, $\lambda=\max(|\lambda_2|,|\lambda_n|) <1$ and 
		\[\|\mW\boldsymbol{\omega}-\mathbf{1}\overline{\omega}\|\le \lambda\|\boldsymbol{\omega}-\mathbf{1}\overline{\omega}\|\]
		for all $\boldsymbol{\omega}\in\mathbb{R}^{n\times d}$, where $\overline{\omega} = \frac{1}{n}\mathbf{1}^{\T}\boldsymbol{\omega}$.
	\end{lemma}
	As a result, when running a consensus algorithm (which is just (\ref{eq: z_i,k}) without gradient descent)
	\begin{equation}
		\label{consensus}
		z_i(k+1) = \sum_{j=1}^{n}w_{ij}z_j(k),
	\end{equation}
	the speed of reaching consensus is determined by $\lambda = \max(|\lambda_2|, |\lambda_n|)$. In particular, if we adopt the so-called lazy Metropolis rule for defining the weights, the dependency of $\lambda$ on the network size $n$ is upper bounded by $1-c/n^2$ for some constant $c$ \cite{olshevsky2017linear}.
	\vspace{1em}
	
	\fbox{\begin{minipage}[c]{0.9\textwidth}
			Lazy Metropolis rule for constructing $\mW$:
			\begin{equation*}
				w_{ij}=\begin{cases}
					\frac{1}{2\max\{\degree(i),\degree(j)\}}, & \text{if }i\in \mathcal{N}_i,  \\
					1- \sum_{j\in\mathcal{N}_i}w_{ij}, & \text{if }i=j,\\
					0, & \text{otherwise}.
				\end{cases}
			\end{equation*}
			Notation: $\degree(i)$ denotes the degree (number of ``neighbors'') of node $i$.
			Correspondingly, $\mathcal{N}_i$ is the set of ``neighbors'' for agent $i$.
	\end{minipage}}
	\vspace{1em}
	
	Despite the fact that $\lambda$ may be very close to $1$ with large $n$, the consensus algorithm (\ref{consensus}) enjoys geometric convergence speed, i.e., 
	\[\sum_{i=1}^n\left\|z_i(k)-\frac{1}{n}\sum_{j=1}^n z_j(k)\right\|^2 \le \lambda^k\sum_{i=1}^n\left\|z_i(0)-\frac{1}{n}\sum_{j=1}^n z_j(0)\right\|^2.\] 
	By contrast, the optimal rate of convergence for any stochastic gradient methods is sublinear, asymptotically $\mathcal{O}(\frac{1}{k})$ (see \cite{nemirovski2009robust}). This difference suggests that a consensus-based distributed algorithm for stochastic optimization may match the centralized methods in the long term: any errors due to consensus will decay at a fast-enough rate so that they ultimately do not matter. 
	
	In what follows, we discuss and compare the performance of the centralized stochastic gradient descent (SGD) method and DSGD. We will show that both methods asymptotically converge at the rate $\frac{\sigma^2}{n\mu^2 k}$. Furthermore, the time needed for DSGD to approach the asymptotic convergence rate turns out to scale as  $\mathcal{O} \left(\frac{n}{(1-\lambda)^2} \right)$.
	
	\subsection{Centralized Stochastic Gradient Descent (SGD)}
	
	The benchmark for evaluating the performance of DSGD is the centralized stochastic gradient descent (SGD) method, which we now  describe. At each iteration $k$, the following update is executed:
	\begin{equation}
		\label{eq: centralized}
		z(k+1)=z(k)-\alpha_k \bar{g}(k),
	\end{equation}
	where stepsizes satisfy $\alpha_k=\frac{1}{\mu k}$,
	and $\bar{g}(k)=\frac{1}{n}\sum_{i=1}^n g_i(z(k),\xi_i(k))$,
	i.e., $\bar{g}(k)$ is the average of $n$ noisy gradients evaluated at $z(k)$ (by utilizing $n$ gradients at each iteration, we are keeping the computational power the same for  SGD and DSGD). As a result, the gradient estimation is more accurate than using just one gradient. Indeed, from Assumption \ref{asp: gradient samples} we have
	\begin{equation}
		\label{SGD_gradient_samples}
		\bE[\|\bar{g}(k)-\nabla f(z(k))\|^2]=\frac{1}{n^2}\sum_{i=1}^n\bE[\|g_i(z(k),\xi_i(k))-\nabla f_i(z(k))\|^2]\le \frac{\sigma^2}{n}.
	\end{equation}
	
	We measure the performance of SGD by $R(k)=\bE[\|z(k)-z^*\|^2]$, the expected squared distance between the solution at time $k$ and the optimal solution.
	Theorem \ref{Thm: centralized} characterizes the convergence rate of $R(k)$, which is optimal for such stochastic gradient methods (see \cite{nemirovski2009robust,rakhlin2012making}).
	\begin{theorem}
		\label{Thm: centralized}
		Under SGD (\ref{eq: centralized}), supposing Assumptions \ref{asp: mu-L_convexity}-\ref{asp: network} hold, we have
		\begin{align}
			\label{centralized_rate}
			R(k)
			\leq  \frac{\sigma^2}{n\mu^2 k}+\mathcal{O}_k\left(\frac{1}{k^2}\right).
		\end{align}
	\end{theorem}
	To compare with the analysis for DSGD later, we briefly describe how to obtain (\ref{centralized_rate}). Note that 
	\begin{equation*}
		R(k+1)= \bE[\|z(k)-\alpha_k \bar{g}(k)-z^*\|^2]
		= \bE[\|z(k)-\alpha_k \nabla f(z(k))-z^*\|^2]+\alpha_k^2\bE[\|\nabla f(z(k))-\bar{g}(k)\|^2].
	\end{equation*}
	For large $k$, in light of Lemma \ref{lem: contraction_mu-L_convexity} and relation (\ref{SGD_gradient_samples}), we have the following inequality that relates $R(k+1)$ to $R(k)$.
	\begin{equation}
		\label{centralized_first_inequality}
		R(k+1) \le (1-\alpha_k\mu)^2 R(k)+\frac{\alpha_k^2\sigma^2}{n}
		=  \left(1-\frac{1}{k}\right)^2 R(k)+\frac{\sigma^2}{n\mu^2}\frac{1}{k^2}.
	\end{equation} 
	A simple induction then gives Eq. (\ref{centralized_rate}). 
	
	\subsection{Distributed Stochastic Gradient Descent (DSGD)}
	\label{subsec: DSGD}
	
	We assume the same stepsize policy for DSGD and SGD.
	To analyze DSGD starting from Eq. (\ref{eq: z_i,k}), define
	\begin{equation}
		\oz(k)=\frac{1}{n}\sum_{i=1}^n z_i(k)
	\end{equation}
	as the average of all the iterates in the network.
	Differently from the analysis for SGD, we will be concerned with two error terms. The first term $\bE[\|\oz(k)-z^*\|^2]$, called the expected optimization error, defines the expected squared distance between $\oz(k)$ and $z^*$, and the second term $\sum_{i=1}^n\bE\left[\|z_i(k)-\oz(k)\|^2\right]$, called the expected consensus error, measures the dissimilarities of individual estimates among all the agents.  
	The average squared distance between individual iterate $z_i(k)$ and the optimum $z^*$ is given by
	\begin{equation}
		\label{Error_decouple}
		\frac{1}{n}\sum_{i=1}^n\bE\left[\|z_i(k)-z^*\|^2\right]=\bE[\|\oz(k)-z^*\|^2]+\frac{1}{n}\sum_{i=1}^n\bE\left[\|z_i(k)-\oz(k)\|^2\right].
	\end{equation}
	\noindent 
	Hence, exploring the two terms will provide us with insights into the performance of DSGD.
	To simplify notation, denote $U(k)=\bE\left[\|\oz(k)-z^*\|^2\right]$, $V(k)=\sum_{i=1}^n\bE\left[\|z_i(k)-\oz(k)\|^2\right]$, $\forall k$.
	
	Inspired by the analysis for SGD, we first look for an inequality that bounds $U(k)$, which is analogous to $\bE[\|z(k)-z^*\|^2]$ in SGD. One such relation turns out to be \cite{pu2019non}:
	\begin{equation}
		\label{Opt_error_pre}
		U(k+1)	\le \left(1-\frac{1}{k}\right)^2U(k)
		+\frac{2 L}{\sqrt{n}\mu}\frac{\sqrt{U(k)V(k)}}{k}
		+\frac{L^2}{n\mu^2}\frac{V(k)}{k^2}+\frac{\sigma^2}{n\mu^2}\frac{1}{k^2}.
	\end{equation}
	Comparing (\ref{Opt_error_pre}) to (\ref{centralized_first_inequality}), we find two additional terms on the right-hand side of the inequality. Both terms involve the expected consensus error $V(k)$, thus reflecting the additional disturbances caused by the dissimilarities of solutions. Relation (\ref{Opt_error_pre}) also suggests that the convergence rate of $U(k)$ can not be better than $R(k)$ for SGD, which is expected. Nevertheless,  if $V(k)$ decays fast enough compared to $U(k)$, it is likely that the two additional terms are negligible in the long run, and we would guess  that the convergence rate of $U(k)$ is comparable to $R(k)$ for SGD. 
	
	This indeed turns out to be the case, as it is shown in \cite{pu2019non} that $V(k) \le \mathcal{O} \left(\frac{n}{(1-\lambda)^2}\right)\frac{1}{k^2}$ for $k\ge \mathcal{O} \left(\frac{1}{(1-\lambda)}\right)$. Plugging this into Eq. (\ref{Opt_error_pre}) leads to the inequality $U(k) \le  \frac{\theta^2\sigma^2}{(1.5\theta-1)n\mu^2 k}+ \mathcal{O} \left(\frac{1}{(1-\lambda)^2}\right)\frac{1}{k^2}$.
	Hence, when $k \geq \mathcal{O} \left(\frac{n}{(1-\lambda)^2}\right)$,
	we have that 
	\[\frac{1}{n}\sum_{i=1}^n\bE\left[\|z_i(k)-z^*\|^2\right] \leq \frac{\sigma^2}{n \mu^2 k} \mathcal{O}(1). \] In other words, we have the asymptotic network independence phenomenon: after a transient, DSGD performs comparably to a centralized stochastic gradient descent method with the same computational power (e.g., which can query the same number of gradients per step as the entire network).

	\subsection{Numerical Illustration}
	
	We provide a numerical example to illustrate the asymptotic network independence property of DSGD. 
	Consider the \emph{on-line} Ridge regression problem
	\begin{equation}
		\label{Ridge Regression}
		z^*=\arg\min_{z\in\mathbb{R}^d}\sum_{i=1}^n f_i(z)\left(=\mathbb{E}_{u_i,v_i}\left[\left(u_i^{\T} z-v_i\right)^2+\rho\|z\|^2\right]\right),
	\end{equation}
	where 
	$\rho>0$ is a penalty parameter.
	Each agent $i$ collects data points in the form of $(u_i,v_i)$ continuously over time with $u_i\in\mathbb{R}^d$ representing the features and $v_i\in\mathbb{R}$ being the observed outputs. Suppose each $u_i\in[-1,1]^d$ is uniformly distributed, and $v_i$ is drawn according to $v_i=u_i^{\T} \tilde{z}_i+\varepsilon_i$, where $\tilde{z}_i$ are predefined parameters uniformly situated in $[0,10]^d$, and $\varepsilon_i$ are independent Gaussian random variables with mean $0$ and variance $1$.
	Given a pair $(u_i,v_i)$, agent $i$ can compute an estimated gradient of $f_i(z)$: $g_i(z,u_i,v_i)=2(u_i^{\T}z -v_i)u_i+2\rho z$, which is unbiased.
	Problem (\ref{Ridge Regression}) has a unique solution $z^*$ given by 
	$z^*=\left(\sum_{i=1}^n\mathbb{E}_{u_i}[u_i u_i^{\T}]+n\rho\mathbf{I}\right)^{-1}\sum_{i=1}^n\mathbb{E}_{u_i}[u_i u_i^{\T}]\tilde{z}_i$.
	
	In the experiments, we consider two instances. In the first instance, we assume $n=50$ agents constitute a random network for DSGD, where every two agents are linked with probability $0.2$. In the second instance,  we let $n=49$ agents form a $7\times 7$ grid network. We use Metropolis weights in both instances. 
	The problem dimension is set to $d=10$ and $z_i(0)=\mathbf{0}$, the zero vector, for all $i$. The penalty parameter is set to $\rho=0.1$ and the stepsizes $\alpha_k=\frac{5}{k}$.  
	For both SGD and DSGD, we run the simulations $100$ times and average the results to approximate the expected errors.
	\begin{figure}[htbp]
		\centering
		\begin{subfigure}{0.49\textwidth}
			\includegraphics[width=\textwidth]{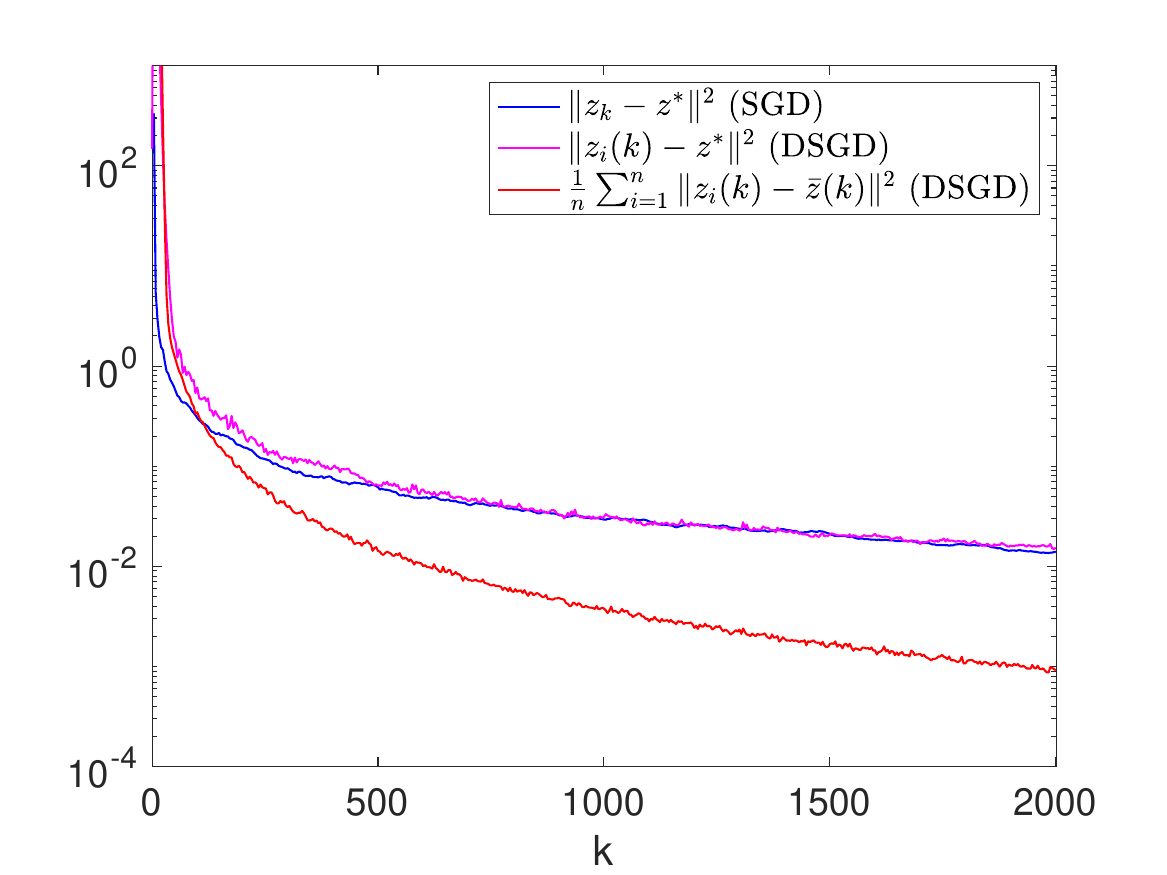}
			\caption{Instance 1 (random network for DSGD).}
		\end{subfigure}
		\begin{subfigure}{0.49\textwidth}
			\includegraphics[width=\textwidth]{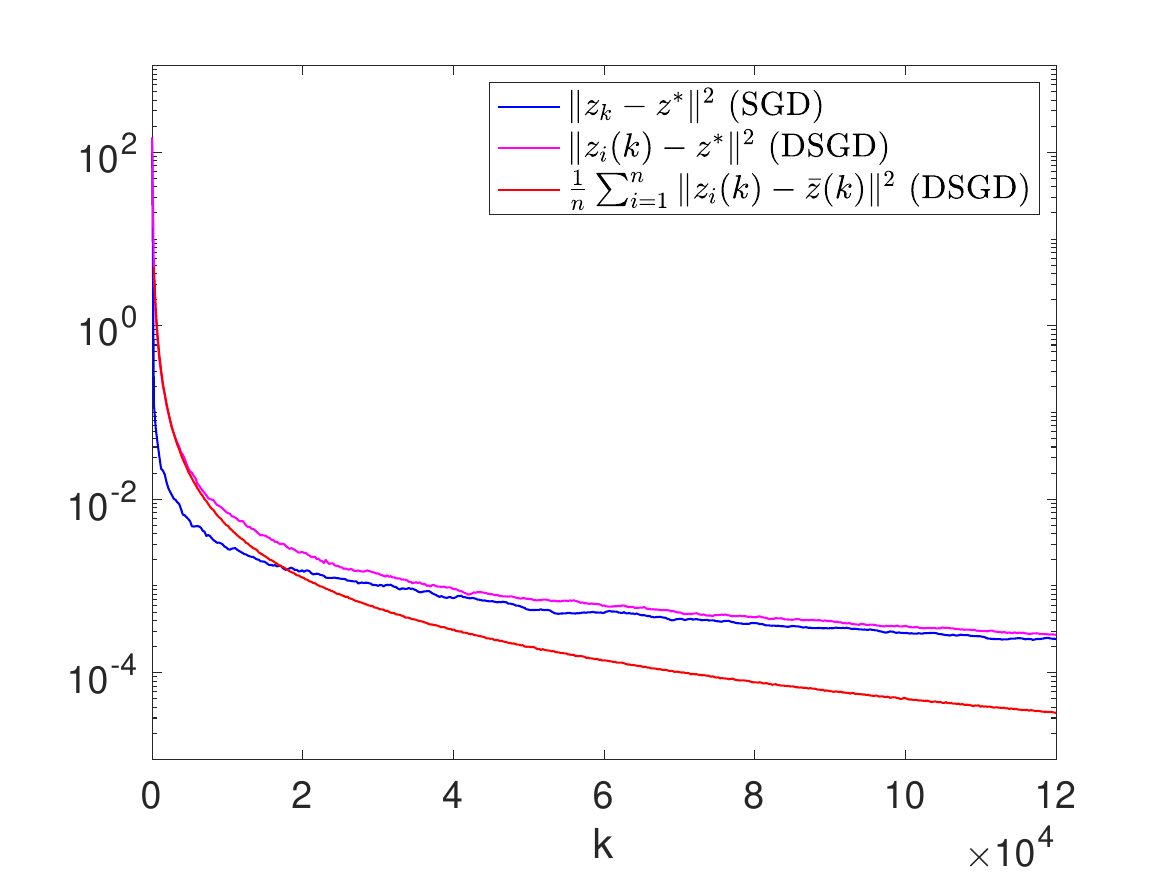}
			\caption{Instance 2 (grid network for DSGD).}
		\end{subfigure}
		\caption{Performance comparison between DSGD and SGD for on-line Ridge regression. For DSGD, the plots show the iterates generated 
			by a randomly selected node $i$ from the set $\{1,2,\ldots,n\}$. Results are averaged over $100$ Monte-Carlo simulations.}
		\label{fig: comparison}
	\end{figure}
	
	The performance of SGD and DSGD is shown in Figure \ref{fig: comparison}. We notice that in both instances the expected consensus error for DSGD converges to $0$ faster than the expected optimization error as predicted from our previous discussion. Regarding the expected optimization error, DSGD is slower than SGD in the first $\sim 800$ (resp., $\sim 4\times 10^4$) iterations for random network (resp., grid network). 
	But after that, their performance is almost indistinguishable. The difference in the transient times is due to the stronger connectivity (or smaller $\lambda$) of the random network compared to the grid network.
	
	\section{Conclusions}
	
	
	In this paper, we provided a discussion of recent results which have overcome a
	barrier in distributed stochastic optimization methods for machine learning under certain scenarios. These results established an asymptotic
	network independence property, that is, asymptotically, the distributed algorithm achieves comparable convergence rate to a centralized algorithm
	with the same computational power.
	We explain the property by
	examples of training ML models and provide a short mathematical analysis.
	
	Along the line of achieving asymptotic network independence in distributed optimization, there are various future research directions, including considering nonconvex objective functions, reducing communication costs and transient time, and using exact gradient information. We next briefly discuss these. 
	
	First, distributed training of deep neural networks --- the state-of-the-art machine learning approach in many application areas --- involves minimizing nonconvex objective functions which are different from the main objectives considered in this paper. This area is largely unexplored with a few recent works in \cite{morral2017success,lian2017can,scaman2019optimal,assran2018stochastic}.
	
	In distributed algorithms, the costs associated with communication among the agents are often non-negligible and may become the main burden for large networks. It is therefore important to explore communication reduction techniques that do not sacrifice the asymptotic  network independence property. The recent papers \cite{assran2018stochastic,koloskova2019decentralized}  touched upon this point.
	
	When considering asymptotic network independence for distributed optimization, an important factor is the transient time to reach the asymptotic convergence rate, as it may take a long time before the distributed implementation catches up with the corresponding centralized method. In fact, as we have shown in Section \ref{subsec: setup}, this transient time can be a function of the network topology and grows with the network size. Reducing the transient time is thus a key future objective.
	
	Finally, while several recent works have established the asymptotic network independence property in distributed optimization, they are mainly constrained to using stochastic gradient information. If the exact gradient is available, can distributed methods   compete with the centralized ones?
	As we know, centralized algorithms typically enjoy a faster convergence speed with exact gradients. For example, plain gradient descent achieves linear convergence for strongly convex and smooth objective functions. To the best of the authors' knowledge, as of writing this survey, with the exception of \cite{li2019decentralized, scaman2019optimal}, results on asymptotic network independence in this setting are currently lacking.

	%
	%
	
	\section*{Acknowledgments}
	We would like to thank Artin Spiridonoff from Boston University for his kind help in providing Figure \ref{SVM}.
	
	{\small\bibliography{ieeesp}{}

\begin{thebibliography}{10}

\bibitem{assran2018stochastic}
M.~Assran, N.~Loizou, N.~Ballas, and M.~Rabbat.
\newblock Stochastic gradient push for distributed deep learning.
\newblock In {\em International Conference on Machine Learning}, pages
  344--353, 2019.

\bibitem{brisimi2018federated}
T.~S. Brisimi, R.~Chen, T.~Mela, A.~Olshevsky, I.~C. Paschalidis, and W.~Shi.
\newblock Federated learning of predictive models from federated electronic
  health records.
\newblock {\em International Journal of Medical Informatics}, 112:59--67, 2018.

\bibitem{chen2012diffusion}
J.~Chen and A.~H. Sayed.
\newblock Diffusion adaptation strategies for distributed optimization and
  learning over networks.
\newblock {\em IEEE Transactions on Signal Processing}, 60(8):4289--4305, 2012.

\bibitem{chen2012limiting}
J.~Chen and A.~H. Sayed.
\newblock On the limiting behavior of distributed optimization strategies.
\newblock In {\em 2012 50th Annual Allerton Conference on Communication,
  Control, and Computing (Allerton)}, pages 1535--1542. IEEE, 2012.

\bibitem{chen2015learning}
J.~Chen and A.~H. Sayed.
\newblock On the learning behavior of adaptive networks—part i: Transient
  analysis.
\newblock {\em IEEE Transactions on Information Theory}, 61(6):3487--3517,
  2015.

\bibitem{chen2015learning2}
J.~Chen and A.~H. Sayed.
\newblock On the learning behavior of adaptive networks—part ii: Performance
  analysis.
\newblock {\em IEEE Transactions on Information Theory}, 61(6):3518--3548,
  2015.

\bibitem{chen2018robust}
R.~Chen and I.~C. Paschalidis.
\newblock A robust learning approach for regression models based on
  distributionally robust optimization.
\newblock {\em The Journal of Machine Learning Research}, 19(1):517--564, 2018.

\bibitem{durrett2007random}
R.~Durrett.
\newblock {\em Random Graph Dynamics}, volume 200.
\newblock Cambridge university press Cambridge, 2007.

\bibitem{koloskova2019decentralized}
A.~Koloskova, S.~U. Stich, and M.~Jaggi.
\newblock Decentralized stochastic optimization and gossip algorithms with
  compressed communication.
\newblock {\em Proceedings of Machine Learning Research}, 97(CONF), 2019.

\bibitem{lan2017communication}
G.~Lan, S.~Lee, and Y.~Zhou.
\newblock Communication-efficient algorithms for decentralized and stochastic
  optimization.
\newblock {\em Mathematical Programming}, pages 1--48, 2017.

\bibitem{li2019decentralized}
Z.~Li, W.~Shi, and M.~Yan.
\newblock A decentralized proximal-gradient method with network independent
  step-sizes and separated convergence rates.
\newblock {\em IEEE Transactions on Signal Processing}, 67(17):4494--4506,
  2019.

\bibitem{lian2017can}
X.~Lian, C.~Zhang, H.~Zhang, C.-J. Hsieh, W.~Zhang, and J.~Liu.
\newblock Can decentralized algorithms outperform centralized algorithms? a
  case study for decentralized parallel stochastic gradient descent.
\newblock In {\em Advances in Neural Information Processing Systems}, pages
  5336--5346, 2017.

\bibitem{morral2014success}
G.~Morral, P.~Bianchi, and G.~Fort.
\newblock Success and failure of adaptation-diffusion algorithms for consensus
  in multi-agent networks.
\newblock In {\em 53rd IEEE Conference on Decision and Control}, pages
  1476--1481. IEEE, 2014.

\bibitem{morral2017success}
G.~Morral, P.~Bianchi, and G.~Fort.
\newblock Success and failure of adaptation-diffusion algorithms with decaying
  step size in multiagent networks.
\newblock {\em IEEE Transactions on Signal Processing}, 65(11):2798--2813,
  2017.

\bibitem{nedic2016stochastic}
A.~Nedi{\'c} and A.~Olshevsky.
\newblock Stochastic gradient-push for strongly convex functions on
  time-varying directed graphs.
\newblock {\em IEEE Transactions on Automatic Control}, 61(12):3936--3947,
  2016.

\bibitem{nedic2009distributed}
A.~Nedic, A.~Olshevsky, A.~Ozdaglar, and J.~N. Tsitsiklis.
\newblock On distributed averaging algorithms and quantization effects.
\newblock {\em IEEE Transactions on Automatic Control}, 54(11):2506--2517,
  2009.

\bibitem{nor}
A.~Nedi{\'c}, A.~Olshevsky, and M.~Rabbat.
\newblock Network topology and communication-computation tradeoffs in
  decentralized optimization.
\newblock {\em Proceedings of the IEEE}, 106(5):953--976, 2018.

\bibitem{nedic2009distributed2}
A.~Nedic and A.~Ozdaglar.
\newblock Distributed subgradient methods for multi-agent optimization.
\newblock {\em IEEE Transactions on Automatic Control}, 54(1):48--61, 2009.

\bibitem{nemirovski2009robust}
A.~Nemirovski, A.~Juditsky, G.~Lan, and A.~Shapiro.
\newblock Robust stochastic approximation approach to stochastic programming.
\newblock {\em SIAM Journal on Optimization}, 19(4):1574--1609, 2009.

\bibitem{olshevsky2017linear}
A.~Olshevsky.
\newblock Linear time average consensus and distributed optimization on fixed
  graphs.
\newblock {\em SIAM Journal on Control and Optimization}, 55(6):3990--4014,
  2017.

\bibitem{olshevsky2018robust}
A.~Olshevsky, I.~C. Paschalidis, and A.~Spiridonoff.
\newblock Robust asynchronous stochastic gradient-push: asymptotically optimal
  and network-independent performance for strongly convex functions.
\newblock {\em arXiv preprint arXiv:1811.03982}, 2018.

\bibitem{pu2017flocking}
S.~Pu and A.~Garcia.
\newblock A flocking-based approach for distributed stochastic optimization.
\newblock {\em Operations Research}, 66(1):267--281, 2017.

\bibitem{pu2018swarming}
S.~Pu and A.~Garcia.
\newblock Swarming for faster convergence in stochastic optimization.
\newblock {\em SIAM Journal on Control and Optimization}, 56(4):2997--3020,
  2018.

\bibitem{pu2018distributed}
S.~Pu and A.~Nedi{\'c}.
\newblock Distributed stochastic gradient tracking methods.
\newblock {\em arXiv preprint arXiv:1805.11454}, 2018.

\bibitem{pu2019non}
S.~Pu, A.~Olshevsky, and I.~C. Paschalidis.
\newblock A non-asymptotic analysis of network independence for distributed
  stochastic gradient descent.
\newblock {\em arXiv preprint arXiv:1906.02702}, 2019.

\bibitem{qu2017harnessing}
G.~Qu and N.~Li.
\newblock Harnessing smoothness to accelerate distributed optimization.
\newblock {\em IEEE Transactions on Control of Network Systems}, 2017.

\bibitem{rakhlin2012making}
A.~Rakhlin, O.~Shamir, and K.~Sridharan.
\newblock Making gradient descent optimal for strongly convex stochastic
  optimization.
\newblock In {\em Proceedings of the 29th International Coference on
  International Conference on Machine Learning}, pages 1571--1578. Omnipress,
  2012.

\bibitem{sayin2017stochastic}
M.~O. Sayin, N.~D. Vanli, S.~S. Kozat, and T.~Ba{\c{s}}ar.
\newblock Stochastic subgradient algorithms for strongly convex optimization
  over distributed networks.
\newblock {\em IEEE Transactions on Network Science and Engineering},
  4(4):248--260, 2017.

\bibitem{scaman2019optimal}
K.~Scaman, F.~Bach, S.~Bubeck, Y.~T. Lee, and L.~Massouli{\'e}.
\newblock Optimal convergence rates for convex distributed optimization in
  networks.
\newblock {\em Journal of Machine Learning Research}, 20(159):1--31, 2019.

\bibitem{sirb2018decentralized}
B.~Sirb and X.~Ye.
\newblock Decentralized consensus algorithm with delayed and stochastic
  gradients.
\newblock {\em SIAM Journal on Optimization}, 28(2):1232--1254, 2018.

\bibitem{towfic2016excess}
Z.~J. Towfic, J.~Chen, and A.~H. Sayed.
\newblock Excess-risk of distributed stochastic learners.
\newblock {\em IEEE Transactions on Information Theory}, 62(10):5753--5785,
  2016.

\bibitem{xin2019variance}
R.~Xin, U.~A. Khan, and S.~Kar.
\newblock Variance-reduced decentralized stochastic optimization with gradient
  tracking.
\newblock {\em arXiv preprint arXiv:1909.11774}, 2019.

\bibitem{ying2018supervised}
B.~Ying, K.~Yuan, and A.~H. Sayed.
\newblock Supervised learning under distributed features.
\newblock {\em IEEE Transactions on Signal Processing}, 67(4):977--992, 2018.

\end{thebibliography}
		\bibliographystyle{abbrv}}
	
	
\end{document}